\documentclass{amsart}
\usepackage[german,english]{babel}
\usepackage[utf8]{inputenc}
\usepackage{newlfont}
\usepackage{amsmath}
\usepackage{amscd}
\usepackage{amssymb}
\usepackage{amsthm}
\usepackage{mathrsfs}
\usepackage{graphicx}
\usepackage{stmaryrd}
\usepackage{amsmath}
\usepackage{latexsym}
\usepackage{braket}
\usepackage{enumerate}
\usepackage{booktabs}
\usepackage{array}
\usepackage{paralist}
\usepackage{pgf,tikz}
\usetikzlibrary{arrows}
\usetikzlibrary{matrix}
\usetikzlibrary{calc}

\usepackage[a4paper,top=3cm,bottom=3cm,left=3.5cm,right=3.5cm,%
bindingoffset=5mm]{geometry}

\newtheorem{theorem}{Theorem}[section]
\newtheorem{intheorem}{Theorem}
\newtheorem{lemma}[theorem]{Lemma}
\newtheorem{corollary}[theorem]{Corollary}
\newtheorem{proposition}[theorem]{Proposition}

\newtheorem{problem}[theorem]{Problem}
\theoremstyle{definition}
\newtheorem{definition}[theorem]{Definition}
\newtheorem{indefinition}{Definition}
\newtheorem{example}[theorem]{Example}
\newtheorem{inproblem}{Problem}

\theoremstyle{remark}
\newtheorem{remark}[theorem]{Remark}

\def\rem{\begin{remark}}
\def\reme{\end{remark}}
\def\ex{\begin{example}}
\def\exe{\end{example}}
\def\tr{\begin{theorem}}
\def\tre{\end{theorem}}
\def\prop{\begin{proposition}}
\def\prope{\end{proposition}}
\def\df{\begin{definition}}
\def\dfe{\end{definition}}
\def\cor{\begin{corollary}}
\def\core{\end{corollary}}
\def\lem{\begin{lemma}}
\def\leme{\end{lemma}}
\def\pr{\begin{proof}}
\def\pre{\end{proof}}
\def\prob{\begin{problem}}
\def\probe{\end{problem}}

\def\it{\begin{itemize}}
\def\ite{\end{itemize}}

\def\intr{\begin{intheorem}}
\def\intre{\end{intheorem}}

\def\inprob{\begin{inproblem}}
\def\inprobe{\end{inproblem}}

\def\indf{\begin{indefinition}}
\def\indfe{\end{indefinition}}

\DeclareMathOperator{\hkf}{HKF}

\newcommand{\diff}{\operatorname{diff}}

\newcommand{\spec}{\operatorname{Spec}}
\newcommand{\rad}{\operatorname{rad}}
\newcommand{\nil}{\operatorname{nil}}
\newcommand{\im}{\operatorname{im}}

\newcommand{\tf}{\operatorname{tf}}
\newcommand{\red}{\operatorname{red}}
\newcommand{\ann}{\operatorname{Ann}}

\newcommand{\supp}{\operatorname{supp}}
\newcommand{\trian}{\operatorname{\triangle}}
\newcommand{\len}{\operatorname{length}}

\def\t{\mathbb T}

\def\n{\mathbb N}

\def\z{\mathbb Z}

\DeclareRobustCommand\longtwoheadrightarrow {\relbar\joinrel\twoheadrightarrow}
\newcommand*{\longhookrightarrow}{\ensuremath{\lhook\joinrel\relbar\joinrel\rightarrow}}
\DeclareFontFamily{U}  {MnSymbolC}{}
\DeclareSymbolFont{MnSyC}         {U}  {MnSymbolC}{m}{n}
\SetSymbolFont{MnSyC}       {bold}{U}  {MnSymbolC}{b}{n}
\DeclareFontShape{U}{MnSymbolC}{m}{n}{
<-6> MnSymbolC5
<6-7> MnSymbolC6
<7-8> MnSymbolC7
<8-9> MnSymbolC8
<9-10> MnSymbolC9
<10-12> MnSymbolC10
<12-> MnSymbolC12}{}
\DeclareFontShape{U}{MnSymbolC}{b}{n}{
<-6> MnSymbolC-Bold5
<6-7> MnSymbolC-Bold6
<7-8> MnSymbolC-Bold7
<8-9> MnSymbolC-Bold8
<9-10> MnSymbolC-Bold9
<10-12> MnSymbolC-Bold10
<12-> MnSymbolC-Bold12}{}

\DeclareMathSymbol{\cupdot}{\mathbin}{MnSyC}{60}

\DeclareFontFamily{U} {MnSymbolF}{}
\DeclareSymbolFont{mnsymbols} {U} {MnSymbolF}{m}{n}
\DeclareFontShape{U}{MnSymbolF}{m}{n}{
<-6> MnSymbolF5
<6-7> MnSymbolF6
<7-8> MnSymbolF7
<8-9> MnSymbolF8
<9-10> MnSymbolF9
<10-12> MnSymbolF10
<12-> MnSymbolF12}{}
\DeclareFontShape{U}{MnSymbolF}{b}{n}{
<-6> MnSymbolF-Bold5
<6-7> MnSymbolF-Bold6
<7-8> MnSymbolF-Bold7
<8-9> MnSymbolF-Bold8
<9-10> MnSymbolF-Bold9
<10-12> MnSymbolF-Bold10
<12-> MnSymbolF-Bold12}{}

\DeclareMathSymbol{\bigcupdot}{\mathop}{mnsymbols}{34}

\title{Hilbert-Kunz multiplicity of binoids}

\author{Bayarjargal Batsukh and  Holger Brenner}

\address{{\small Bayarjargal Batsukh, Department of Mathematics, School of Art and Sciences, National University of Mongolia, Mongolia}}
\email{{\small bayarjargal@smcs.num.edu.mn}}

\address{{\small Holger Brenner, Institut f\"ur Mathematik, Universit\"at Osnabr\"uck, Albrechtstrasse 28a, 49076 Osnabr\"uck, Germany}}
\email{{\small holger.brenner@uni-osnabrueck.de}}

\begin{document}

\thanks{\textit{Mathematics Subject Classification (2010)}: 13D40, 13F55, 05E40, 14M25.
\\ \indent  \textit{Keywords}: Hilbert-Kunz multiplicity, Binoid}

\begin{abstract}
We prove  in a broad combinatorial setting, namely for finitely generated semipositive cancellative  reduced binoids, that the Hilbert-Kunz multiplicity is a rational number independent of the characteristic.
\end{abstract}

\maketitle

\section*{Introduction}

Let $(R,\mathfrak{m})$ be a commutative Noetherian local ring of dimension $d$ containing a field $K$ of positive characteristic $p$.
For an ideal $I$ and a prime power $q=p^e$ we define the ideal $I^{[q]}=\langle a^q|a\in I\rangle$ which is the ideal generated by the $q$th power of the elements of $I$.
Let $I$ be an $\mathfrak{m}$-primary ideal of $R$ and $M$ a finite $R$-module. Then the $R$-modules $M/I^{[q]}M$ have finite length. The \emph{Hilbert-Kunz function} of $M$ with respect to $I$ is
\[ HKF(I,M)(q)=\len(M/I^{[q]}M). \]
If $M=R,I=\mathfrak{m}$ then we have the classical Hilbert-Kunz function $HKF(\mathfrak{m},R)(q)=HK_R(q)$, introduced by Kunz \cite{Kunz}. He showed that $R$ is regular if and only if $HK_R(q)=q^d$ for all $q$.
In \cite{Monsky}, P. Monsky proved that there is a real constant $c(M)$ such that
\[\len(M/I^{[q]}M)=c(M)q^d+O(q^{d-1}).\]
The \emph{Hilbert-Kunz multiplicity} $e_{HK}(I,M)$ of $M$ with respect to $I$ is
\[e_{HK}(I,M):=\lim_{q\rightarrow\infty}\dfrac{\len(M/I^{[q]}M)}{q^d}. \]
There are many questions related to Hilbert-Kunz function and multiplicity.
\inprob
Is the Hilbert-Kunz multiplicity always a rational number?
\inprobe
\inprob
Is there any interpretation in characteristic 0?
\inprobe
For the following problems the ring comes from a finitely generated $\z$-algebra by reduction modulo $p$. 
\inprob
How does the Hilbert-Kunz multiplicity depend on the characteristic $p$?
\inprobe
\inprob
Does the limit
\[\lim_{p\rightarrow\infty} e_{HK}(I_p,R_p) \] exist?
\inprobe
\inprob[C.Miller]
Does the limit
\[ \lim_{p\rightarrow \infty} \dfrac{\len(R_p/I_p^{[p]})}{p^d}\]
exist?
\inprobe
In most known cases the Hilbert-Kunz multiplicity is a rational number, for example for toric rings (\cite{Watanabe}, \cite{Bruns}), monoid rings (\cite{Eto}),
monomial ideals and binomial hypersurfaces (\cite{Conca}), rings of finite Cohen-Macaulay type (\cite{Seibert}), for invariant rings for the action of a finite group on a polynomial ring (follows from \cite[Theorem 2.7]{Watanabeyoshida}), for
two-dimensional graded rings (\cite{Bre06}, \cite{Trivedi2}).
In \cite{Brenner} it is shown that there exist also irrational Hilbert-Kunz multiplicities.

There are many situations where the Hilbert-Kunz multiplicity is independent of the characteristic $p$.
For example, for toric rings (\cite{Watanabe}) or invariant rings as above the Hilbert-Kunz multiplicity is independent of the characteristic of the base field at
least for almost all prime characteristics. But there are also examples where the Hilbert-Kunz multiplicity depends on the characteristic.

We can ask when the limit of the Hilbert-Kunz multiplicity exists for $p\rightarrow\infty$. If so then this limit is a
candidate for the Hilbert-Kunz multiplicity in characteristic zero. This  leads us to the question of whether a characteristic zero Hilbert-Kunz multiplicity could be defined directly. In all known cases this limit exists.
H. Brenner, J. Li and C. Miller (\cite{BreLiMil}) have observed that in all known cases where
\[\lim_{p\rightarrow\infty}e_{HK}(R_p)\]
exists then this double limit can be replaced by the limit of Problem 5.

If the rings are of a more combinatorial nature, like for example monoid rings (\cite{Eto}, \cite{Watanabe}, \cite{Bruns}), Stanley-Reisner rings and binomial hypersurfaces (\cite{Conca}), we have positive answers for all these problems. 
Also, the proofs in these cases are easier compared to the methods of P. Monsky, C. Han, P. Teixeira (\cite{MonHan}, \cite{MonTei04}, \cite{MonTei06}) or the geometric methods of H. Brenner and V. Trivedi (\cite{Bre06}, \cite{Bre07}, \cite{TriFak03}, \cite{Trivedi1}, \cite{Trivedi2}, \cite{Trivedi07}). 
We want to generalize these results to a broad and unified concept of a combinatorial ring. For that we work with a new combinatorial structure namely binoids (pointed monoid) which were introduced in the thesis of S. Boettger \cite{Simone}.
A binoid $(N,+,0,\infty)$ is a monoid with an absorbing element $\infty$ which means that for every $a\in N$ we have $a+\infty=\infty+a=\infty$.
This concept recovers among others monoid rings and Stanley-Reisner rings. 

In the first four sections we will describe some basic properties of binoids and related objects,
namely $N$-sets, smash products, exact sequences. Also we will give the definition and some properties on the dimension of binoids and their binoid algebras.
In Section 5 we will define the Hilbert-Kunz function and multiplicity of binoids. This function is given by counting the elements in certain residue class binoids,
and not the vector space dimension (or length) of residue class rings. We define the Hilbert-Kunz function and Hilbert-Kunz multiplicity not only for binoids but also for $N_+$-primary ideals of $N$ and a finitely generated $N$-set in the following way:

Let $N$ be a finitely generated, semipositive binoid, $T$ a finitely generated $N$-set and $\mathfrak{n}$ an $N_+$-primary ideal of $N$. Then we call the number
\[\hkf^N(\mathfrak{n},T,q)=\hkf(\mathfrak{n},T,q):=\# T/([q]\mathfrak{n}+T)\]
(where for a finite binoid we do not count $\infty$) the Hilbert-Kunz function of $\mathfrak{n}$ on the $N$-set $T$ at $q$. In particular, for $T=N$ and $\mathfrak{n}=N_+$ we have
\[ \hkf(N,q):=\hkf(N_+,N,q)=\#N/[q]N_+. \]
Note that this function is defined for all $q\in\n$, not only for powers of a fixed prime number. Also note that this residue construction is possible in the category of binoids, not in the category of monoids. The combinatorial Hilbert-Kunz multiplicity is the limit of this function divided by $q^{\dim N}$, provided this limit exists and provided that there is a reasonable notion of dimension.

It turns out that the Hilbert-Kunz function of a binoid for $q=p^e$ is the same as the Hilbert-Kunz function of its binoid algebra over a field of characteristic $p$. 
Hence from here we have a chance to study the above mentioned five problems, by just studying the binoid case. 

Since the Hilbert-Kunz function of binoids is given by just counting elements, not vector space dimensions,
this reveals clearer the combinatorial nature of the problem.
For example we have
\[\dim_K K[N]/K[N_+]^{[q]}=\dim_K K[N/[q]N_+]=\# N/[q]N_+.  \]
So the computation of the Hilbert-Kunz function and multiplicity of a binoid is in some sense easier than the standard Hilbert-Kunz function and multiplicity. 
The following is our main combinatorial theorem (Theorem \ref{f.g,s.p,c,r binoid}).

\intr
Let $N$ be a finitely generated, semipositive, cancellative, reduced binoid and $\mathfrak{n}$ be an $N_+$-primary ideal of $N$.
Then $e_{HK}(\mathfrak{n},N)$ exists and is a rational number.
\intre

Note that this is a characteristic-free statement and that the existence does not follow from Monsky's theorem.
From this result we can deduce a positive answers to our five problems, see Theorem \ref{Miller conjecture} and Theorem \ref{existsrational} in the current setting.

The strategy to prove these theorems is to reduce it step by step to the corresponding results of lower dimensional binoids fulfilling further properties.
The case of primary ideals in a normal toric setting was given by Eto, Bruns, Watanabe and is true without normality. The components of a reduced torsion-free cancellative binoid are toric, and the Hilbert-Kunz multiplicity depends only on the components of maximal dimension.
In the proof of this we encounter also the non-reduced situation, and for this the (strongly) exact sequences of $N$-sets are extremely useful. The reduction from the case
with torsion to the torsion-free case requires a deeper understanding of the torsion-freefication and its relation to the smash product of the torsion-free binoid and a
finite group. In this setting we need to replace positive (local algebras) by semipositive (corresponding to semilocal algebras), so we also have to generalize the Hilbert-Kunz theory on the ring side.

\section{Binoids and their properties}

In this Section we will introduce (commutative) binoids, describe basic properties of binoids and binoid sets and their properties.
For a general introduction to binoids we refer to \cite{Simone}, where most of the basic concepts were developed. 
We focus on material which is relevant for Hilbert-Kunz theory. 

A \emph{binoid} $(N,+,0,\infty)$ is a commutative monoid with an absorbing element $\infty$, given by the property $x +\infty= \infty$.
We write $N^\bullet$ for the set $N\setminus \{\infty\}$. If  $N^\bullet$ is a monoid, then $N$ is called an \emph{integral} binoid. The set of all units form the \emph{unit group} of $N$, denoted by $N^\times$.
The set of all non-units $N\setminus N^\times$ will be denoted by $N_+$. A binoid $N \neq \{ \infty\}$ is called \emph{semipositive}, if $|N^\times|$ is finite
and \emph{positive}, if $N^\times=\{0\}$. From the corresponding concepts in monoid theory it is clear what a homomorphism of binoids (sending $\infty$ to $\infty$),
what the kernel is and what an ideal, a radical, a prime ideal is. Specific for an ideal $I \subseteq N$ in a binoid is that there is a homomorphism $M \rightarrow M/I$
to the \emph{residue class binoid}, where the ideal is sent to $\infty$ and which is injective elsewhere.

For every ideal $I$ and $q \in \n_+$ we will denote by
\[ [q]I :=\langle q a\mid a\in I\rangle \]
the ideal generated by the set $\{qa\mid a\in I\}$, which we call the $q$th \emph{Frobenius sum} of the ideal.

For a binoid homomorphism $\varphi: N\rightarrow M$ and an ideal $I\subseteq N$ we denote by $I+M$ the ideal generated by $\varphi(I)$
and call it the \emph{extended ideal}. Since $[q]:N\rightarrow N, x\mapsto qx$, is a binoid homomorphism, the ideal $[q]I$ can be considered as
the extended ideal under this homomorphism. The extension of ideals commute with Frobenius sums. An ideal is called a \emph{primary ideal} if $\rad(\mathfrak{n})=N_+$. For an
$N_+$-primary ideal $\mathfrak{n}$ also its Frobenius sums $[q]\mathfrak{n}$ are $N_+$-primary. For a finitely generated semipositive binoid and an $N_+$-primary ideal $\mathfrak{n}$
the residue class binoid $N/\mathfrak{n}$ is a finite set.

An element $a\in N$ is called \emph{nilpotent} if $na=a+\cdots+a=\infty$ for some $n\in\n$. The set of all nilpotent elements will be denoted by $\nil(N)$.
We say that $N$ is  \emph{reduced} if $\nil(N)=\{\infty\}$. We call the quotient binoid $N_{\red}:=N/\nil(N)$ the \emph{reduction} of $N$.
An element $a \in N$ is a \emph{torsion} element in case $a=\infty$ or $na = nb$ for some $b\in N, b\neq a, n \geqslant 2$. 
We say $N$ is \emph{torsion-free} if there are no other torsion elements in M besides $\infty$, i.e. $na = nb$ implies $a = b$ for every $a, b \in N$ and $n\geqslant 1$.
A binoid is called \emph{torsion-free up to nilpotence} if $na = nb \neq \infty$ implies $a = b$ for
every $a, b \in N$ and $n \geqslant 1$.

For an integral binoid $N$ we denote by $\diff N$ the \emph{difference} binoid of $N$, which is a group binoid. If $N$ is integral and cancellative then there is an injection
$N \subseteq \diff N$. If $N \subseteq M \subseteq \diff N$ we say that $M$ is \emph{birational} over $N$.

The \emph{spectrum} of $N$, denoted by $\spec N$, is the set of all prime ideals of $N$. It can be made to a (finite, if $N$ is finitely generated binoid) topological space.
The \emph{combinatorial dimension} of a binoid $N$, denoted by $\dim N$, is the supremum of the lengths of strictly increasing chains of prime ideals of $N$. Without any further
condition the Krull dimension of a binoid algebra $K[N]$ over a field $K$ and $\dim N$ need not be the same. A binoid algebra is the monoid algebra where one additionally sets $T^\infty =0$.

\lem
\label{dimension properties}
Let $N$ be a finitely generated binoid. If $N$ is integral and $I\neq \{\infty\}$ be an ideal of $N$, then $\dim N/I< \dim N$. If $\mathfrak{p}$ and $\mathfrak{q}$ are
different minimal prime ideals of $N$, then $\dim N/(\mathfrak{p}\cup \mathfrak{q})<\min\{\dim N/\mathfrak{p},\dim N/\mathfrak{q}\}$.
\leme
\pr
See \cite[Proposition 1.8.3]{Batsukhthesis} and \cite[Proposition 1.8.4]{Batsukhthesis}.
\pre

\section{$N$-sets}

A \emph{pointed set} $(S,p)$ is a set $S$ with a distinguished element $p\in S$. Let $N$ be a binoid.
\df
An \emph{operation} of $N$ on a pointed set $(S,p)$ is a map
\[ +:N\times S \longrightarrow S,~~ (n,s)\longmapsto n+s, \]
such that the following conditions are fulfilled:
\begin{enumerate}
  \item For all $n,m\in N$ and $s\in S:(n+m)+s=n+(m+s)$.
  \item For all $s\in S: 0+s=s$.
  \item For all $s\in S: \infty+s=p$.
  \item For all $n\in N: n+p=p$.
\end{enumerate}
Then $S$ is called an $N$-\emph{set}.
\dfe

Given a fixed binoid homomorphism $N \rightarrow M$ then $M$ is an $N$-set in a natural way. This applies in particular for $N$ itself and for residue class binoids.
For an $N$-set $(S,p)$ we say that $a\in N$ is an \emph{annihilator} if $a+s=p$ for every $s\in S$. We denote the set of all annihilators by $\ann S$, which is an ideal of $N$.
The $N$-set $(S,p)$ is also an $(N/\mathfrak{a})$-set for every ideal $\mathfrak{a}\subseteq \ann S$. If there exist finitely many elements
$s_1,\dots,s_k\in S$ such that for all $s\in S$ we can write $s=n+s_j$ where $n\in N$,
then we say that $S$ is a \emph{finitely generated $N$-set}. We call the elements $s_j$ $N$-\emph{generators}.

For a finite $N$-set $S$ we set
\[\# S=|S|-1=|S\setminus\{p\}| , \]
so we do not count the distinguished point.

For a family of $N$-sets $(S_i,p_i)$, $i\in I$, we define the \emph{pointed union} of $S_i,i\in I,$ by
\[\bigcupdot_{i\in I} S_i=(\biguplus_{i\in I}S_i)/\sim, \]
where $\biguplus$ is the disjoint union and $a\sim b$ if and only if $a=b$ or $a=p_j,~b=p_k$ for some $j,k$.
So the pointed union just contracts the points $p_j$ to one point. We write $S^{\cupdot r}=\bigcupdot_{i=1}^r S$ and in particular
$N^{\cupdot r}$ for the $r$-folded pointed union of $N$ with itself.

For a homomorphism $f:S\rightarrow T$ of $N$-sets we set
\[ \im(f)=\{t\in T \mid t = f(s) \text{ for some } s \in S\} \text{ and } \ker(f)=\{s\in S \mid  f(s)=p_T\}. \]
For an $N$-subset $S\subseteq T$ we define the \emph{quotient} of $T$ by $S$ to be the $N$-set $(T \setminus S) \cup {p}$ and denoted it by $T/S$,
so $S$ is contracted to a point.

The following statements are analogous to statements about modules, for the easy proofs we refer to \cite{Batsukhthesis}.

\lem
\label{hom}
Let $N$ be a binoid, $(S,p_S),(T,p_T)$ be $N$-sets and $S'\subseteq S$ an $N$-subset. If we have an $N$-set homomorphism $\phi:S\rightarrow T$ with $\phi(S')=p_T$ then there exists a unique  homomorphism $\tilde{\phi}:S/S'\rightarrow T$ such that the following diagram commutes.

\[\begin{tikzpicture}[node distance=1.8cm, auto]
  \node (S) {$S$};
  \node (T) [right of=S] {$T$};
  \node (A) [below of=S] {$S/S'$};
  \draw[->] (S) to node {$\phi$} (T);
  \draw[->, dashed] (A) to node  [swap]{$\tilde{\phi}$} (T);
  \draw[->] (S) to node [swap] {$\varphi$} (A);
\end{tikzpicture}
\]
If $\phi$ is surjective, then $\tilde{\phi}$ is surjective.
\leme
\pr
See \cite[Lemma 1.5.7]{Batsukhthesis}.
\pre

\prop
\label{quotient to union}
Let $N$ be a finitely generated binoid, $J$ be an ideal of $N$ and $S\subseteq T$ be $N$-sets. Then $(T/S)/(J+(T/S))= T/(S\cup (J+T))$.
\prope
\pr
See \cite[Proposition 1.5.8]{Batsukhthesis}.
\pre

\lem
\label{canonical isomorphism}
Let $N$ be a binoid, $I\subseteq N$ an ideal of $N$ and $(T,p)$ be an $N$-set.
If $r$ is some positive integer then we have a canonical $N$-set isomorphism
\[T^{\cupdot r}/(I+T^{\cupdot r})\cong (T/(I+T))^{\cupdot r}. \]
\leme
\pr
See \cite[Lemma 1.5.9]{Batsukhthesis}.
\pre

\lem
\label{surj hom}
Let $N$ be a binoid, $(S,p_S),(T,p_T)$ $N$-sets and $I\subseteq N$ be an ideal of $N$.
If we have a surjective $N$-set homomorphism $\phi:S\rightarrow T$ then there exists a canonical surjective $N$-set homomorphism
$\tilde{\phi}:S/(I+S)\longrightarrow T/(I+T)$.
\leme
\pr
See \cite[Lemma 1.5.10]{Batsukhthesis}.
\pre

\lem
\label{surj map for f.g. set}
Let $N$ be a binoid and $T$ be an $N$-set. Then
\begin{enumerate}
\item For $t_1,\dots,t_r\in T$ we can define an $N$-set homomorphism $\phi:N^{\cupdot r}\rightarrow T$.
\item $t_1,\dots,t_r\in T$ is a generating system of $T$ over $N$ if and only if $\phi$ is a surjective homomorphism.
\item $T$ is finitely generated over $N$ if and only if there exists a surjective $N$-set homomorphism $N^{\cupdot r}\rightarrow T$.
\end{enumerate}
\leme
\pr
See \cite[Lemma 1.5.11]{Batsukhthesis}.
\pre

The smash product of binoids and $N$-sets correspond to the tensor product of algebras and modules.

\df
Let $(S_i,p_i)_{i\in I}$ be a finite family of pointed sets and $\sim_\wedge$ the relation on $\prod_{i\in I} S_i$ given by
\[ (s_i)_{i\in I}\sim_\wedge (t_i)_{i\in I}:\Leftrightarrow s_i=t_i,\text{ for all } i\in I, \text{ or } s_j=p_j,t_k=p_k \text{ for  some } j,k\in I. \]
Then the pointed set
\[\bigwedge_{i\in I}S_i:=(\prod_{i\in I} S_i)/\sim_\wedge \]
with distinguished point $[p_\wedge:=(p_i)_{i\in I}]$ is called the \emph{smash product} of the family $S_i,i\in I$.
\dfe
\df
Let $N$ be a binoid, $(S_i,p_i)_{i\in I}$ be a finite family of pointed $N$-sets and $\sim_{\wedge_N}$ the equivalence relation on
$\bigwedge_{i\in I} S_i$ generated by
\[ \cdots\wedge n+s_i\wedge \cdots \wedge s_j \wedge\cdots \sim_{\wedge_N}\cdots\wedge s_i\wedge\cdots\wedge n+s_j \wedge\cdots,\] for all $i,j\in I$ and $n\in N$.
Then \[\bigwedge_{i\in I} {\!_N} S_i:=(\bigwedge_{i\in I} S_i)/\sim_{\wedge_N} \]
is called the
\textbf{smash product} of the family $(S_i)_{i\in I}$ over $N$.
\dfe

\prop
\label{pointed union and smash}
Let $N$ be a binoid and $S$, $T_i,1\leqslant i\leqslant k$ be $N$-sets. Then
\[S\wedge_N (\bigcupdot_{i=1}^k T_i)=\bigcupdot_{i=1}^k (S\wedge_N T_i). \]
\prope
\pr
See \cite[Proposition 1.6.7]{Batsukhthesis}.
\pre

\prop
\label{quotient}
Let $N$ be a binoid and $(S,p)$ an $N$-set. If $I$ is an ideal of $N$ then
\[(N/I)\wedge_N S\cong S/(I+S).\]
\prope
\pr
See \cite[Proposition 1.6.8]{Batsukhthesis}.
\pre

\prop
\label{surj map to smash}
Let $N$ be a binoid and $J$ an ideal of $N$. If $S$ is a finitely generated $N$-set then there is a surjective homomorphism
$(N/J)^{\cupdot r}\rightarrow S\wedge_N N/J$, where $r$ is the number of generators of $S$.
\prope
\pr
See \cite[Proposition 1.6.9]{Batsukhthesis}.
\pre

\cor
\label{quotient over primary ideal and smash}
Let $N$ be a finitely generated, semipositive binoid and $J$ an $N_+$-primary ideal of $N$. If $S$ is a finitely generated $N$-set then $S\wedge_N N/J$ is finite.
\core
\pr
See \cite[Corollary 1.6.10]{Batsukhthesis}.
\pre

\lem
\label{smash of quotients}
Let $M,N$ be binoids and $I\subseteq M,J\subseteq N$ be ideals. Then $(I\wedge N)\cup (M\wedge J)$ is an ideal of $M\wedge N$ and
\[(M\wedge N)/((I\wedge N)\cup (M\wedge J))\cong M/I\wedge N/J. \]
\leme
\pr
See \cite[Lemma 1.6.13]{Batsukhthesis}.
\pre

\prop
\label{number of elements of smash}
Let $N,M$ be finite binoids, then
\[\#(N\wedge M)=\#N\cdot\#M. \]
\prope
\pr
See \cite[Proposition 1.6.2]{Batsukhthesis}.
\pre

\tr
\label{dimension of smash product}
Let $M,N$ be non zero binoids of finite dimension. Then
\[\dim M \wedge N=\dim M+ \dim N. \]
\tre
\pr
See \cite[Theorem 1.8.5]{Batsukhthesis}.
\pre

\section{Exact sequences for $N$-sets}

The concept of (strongly) exact sequences for $N$-sets is crucial to reduce statements on the Hilbert-Kunz function to lower dimensional binoids.

\df
Let $N$ be a binoid. A sequence 
\[S_0 \xrightarrow{\;\phi_1} S_1 \xrightarrow{\;\phi_2} S_2 \xrightarrow{\;\phi_3}\cdots\xrightarrow{\;\phi_n} S_n \]
of $N$-sets and $N$-set homomorphisms is called \emph{exact} if the image of each homomorphism is equal to the kernel of the next:
\[ \im\phi_k = \ker\phi_{k+1}. \]
\dfe
\df
Let $N$ be a binoid. An exact sequence 
\[S_0 \xrightarrow{\;\phi_1} S_1 \xrightarrow{\;\phi_2} S_2 \xrightarrow{\;\phi_3}\cdots\xrightarrow{\;\phi_n} S_n\] 
is called \emph{strongly exact} if $\phi_k$ is injective on $S_{k-1}\setminus \ker \phi_k$ for every $k$.
\dfe
\prop
\label{exact sequence}
Let $N$ be a binoid, $S\subseteq T$ and $U$ be $N$-sets. Then we have an exact sequence of $N$-sets
 \[\infty\xrightarrow{\;\phi_1} \{s\wedge u\mid s\wedge u=\infty \text{ in } T\wedge_N U\} \xrightarrow{\;\phi_2} S\wedge_N U\xrightarrow{\;\phi_3} T\wedge_N U\xrightarrow{\;\phi_4} (T/S)\wedge_N U\xrightarrow{\;\phi_5}\infty. \]
\prope
\pr
We have an exact sequence
$\infty \rightarrow S\hookrightarrow T\rightarrow T/S\rightarrow \infty$ and we can smash this sequence with $U$. Then we obtain a sequence
\[ \infty\xrightarrow{\;\phi_1} \{s\wedge u\mid s\wedge u=\infty \text{ in } T\wedge_N U\} \xrightarrow{\;\phi_2} S\wedge_N U\xrightarrow{\;\phi_3}
T\wedge_N U\xrightarrow{\;\phi_4} (T/S)\wedge_N U\xrightarrow{\;\phi_5}\infty,\]
where $\phi_2$ is the inclusion, $\phi_3(s\wedge u)=s\wedge u\in T\wedge_N U$, $\phi_4(t\wedge u)=[t]\wedge u$ and $\phi_5([t]\wedge u)=\infty$.
We know by definition that
$\im\phi_1 =\{\infty\}= \ker\phi_{2}$, $\im\phi_2 =\{s\wedge u\mid s\wedge u=\infty \text{ in } T\wedge_N U\}= \ker\phi_{3}$, $\im\phi_3 =S\wedge_N U= \ker\phi_{4}$, $\im\phi_4 =T/S\wedge_N U= \ker\phi_{5}$.
So it is an exact sequence.
\pre
\ex
\label{example of not strongly exact}
Let $N=\n^\infty$ and $S$ be an $N$-set with an operation given by $n+s=f^n(s)$, where $f:S \rightarrow S$ is a pointed map. Then we have an exact sequence
\[\infty\longrightarrow N_+\longrightarrow N\longrightarrow N/N_+\longrightarrow\infty \]
and we can smash this sequence with $S$ over $N$.
Then by Proposition \ref{exact sequence}, we get an exact sequence
\[\infty\rightarrow \{n\wedge s\mid n\geqslant 1, n\wedge s=\infty \text{ in } N\wedge_N S\}\rightarrow N_+\wedge_N S
\rightarrow N\wedge_N S\rightarrow N/N_+\wedge_N S\rightarrow\infty. \]
We have also isomorphisms
\[N\wedge_N S\longrightarrow S,\;n\wedge s\longmapsto n+s=f^n(s),\]
\[N_+\wedge_N S\longrightarrow S,\;n\wedge s\longmapsto (n-1)+s=f^{n-1}(s),\]
and from here we get
$\{n\wedge s\mid n\geqslant 1, n\wedge s=\infty \text{ in } N\wedge_N S\}\cong\{t\in S\mid f(t)=p\}=\ker f$ and
$N/N_+\wedge_N S\cong S/\im f$.
So we have an exact sequence 
\[\infty\longrightarrow \ker f\longrightarrow S\xrightarrow{\;f\;} S \longrightarrow S/\im f\longrightarrow\infty.\] 
If $S=\{a,b,p\}$, $f(a)=f(b)=b$ and $f(p)=p$ then this sequence is not strongly exact.
\exe

\prop
\label{strongly exact sequence}
Let $N$ be a binoid, $J\subseteq N$ an ideal and $S\subseteq T$ be $N$-sets. Then we have a strongly exact sequence of $N$-sets
\[\infty \! \xrightarrow{\phi_1}  \{s\wedge [a]\mid s\wedge [a] \! = \! \infty \text{ in } T\wedge_N N/J\}
\! \xrightarrow{\phi_2} S\wedge_N N/J\xrightarrow{\phi_3} T\wedge_N \! N/J  \! \xrightarrow{\phi_4} T/S\wedge_N N/J
\! \xrightarrow{\phi_5} \! \infty \]
which is the same as 
\[ \infty \! \xrightarrow{\phi_1}(S\cap(J+T))/(J+S) \xrightarrow{\phi_2}
S/(J+S)\xrightarrow{\phi_3}T/(J+T)\xrightarrow{\phi_4}(T/S)/(J+T/S))\xrightarrow{\phi_5} \! \infty. \]
If $S=I$ is an ideal of $N$ and $T=N$, then we have the strongly exact sequence
\[\infty\longrightarrow (I\cap J)/(I+J) \longrightarrow I/(I+J) \longrightarrow N/J\longrightarrow (N/I)/(J+N/I)\longrightarrow\infty. \]
\prope
\pr
From Proposition \ref{exact sequence}, when $U=N/J$, we have an exact sequence 
\[ \infty \! \xrightarrow{\phi_1} \{s\wedge [a]\mid s\wedge [a]=\infty \text{ in } T\wedge_N N/J\}\! \xrightarrow{\phi_2}
S\wedge_N N/J \! \xrightarrow{\phi_3} \! T\wedge_N \!  N/J \! \xrightarrow{\phi_4} T/S\wedge_N \! N/J \! \xrightarrow{\phi_5} \! \infty. \]
By Proposition \ref{quotient} we know that $S\wedge_N N/J\cong S/(J+S)$ and $T\wedge_N N/J\cong T/(J+T)$.
Let
\[s\wedge [a]\in\{s\wedge [a]\mid s\wedge [a]=\infty \text{ in } T\wedge_N N/J\}. \]
Then $S/(J+S)\ni[a+s]=\infty\in T/(J+T)$. So we have $a+s\in J+T$, which means that $a+s\in S\cap (J+T)$. Hence we have
\[\{s\wedge [a]\mid s\wedge [a]=\infty \text{ in } T\wedge_N N/J\}\cong(S\cap(J+T))/(J+S).\]
Also, by Proposition \ref{quotient} and Proposition \ref{quotient to union}, we get 
\[ T/S\wedge_N N/J\cong (T/S)/(J+T/S)=T/(S\cup(J+T)). \]
So we can rewrite the previous exact sequence as
\[\infty\xrightarrow{\phi_1}(S\cap(J+T))/(J+S) \xrightarrow{ \phi_2}S/(J+S)\xrightarrow{\phi_3}T/(J+T)\xrightarrow{\phi_4}(T/S)/(J+T/S))\xrightarrow{\phi_5}\infty \]
or, equivalently
\[ \infty\xrightarrow{\phi_1}(S\cap(J+T))/(J+S) \xrightarrow{
\phi_2}S/(J+S)\xrightarrow{\phi_3}T/(J+T)\xrightarrow{\phi_4}T/(S\cup(J+T))\xrightarrow{\phi_5}\infty. \]
Here $\phi_1(\infty)=\infty$, $\phi_2$ is the inclusion (so it is injective), $\phi_3$ is an inclusion (injective) on $S\setminus (J+T)$, $\phi_4$ is surjective and outside of the kernel it is a bijection, $\phi_5([s])=\infty$. So it means that our sequence is a strongly exact sequence.

If $S=I$ and $T=N$ then we get
\[ \infty\xrightarrow{\;\phi_1}I\cap J/(J+I) \xrightarrow{\;\phi_2}I/(J+I)\xrightarrow{\;\phi_3}N/J\xrightarrow{\;\phi_4}(N/I)/(J+N/I)\xrightarrow{\;\phi_5}\infty.\qedhere \]
\pre

\prop
\label{general equation of exact seq}
Let $N$ be a binoid and $\infty \rightarrow S_1 \rightarrow S_2 \rightarrow\cdots\rightarrow S_n\rightarrow \infty$ a strongly exact sequence of finite $N$-sets. Then
 \[ \sum_{i=1}^n (-1)^i\# S_i=0. \]
\prope
\pr
Write $S_i=K_i\uplus R_i\uplus \{p_i\}$, with maps
 \[ \begin{aligned}
  \phi_i:S_{i-1}&\longrightarrow S_i,\\
      R_{i-1}&\xrightarrow{bij} K_i,\\
      K_{i-1}&\longrightarrow p_i,\\
      p_{i-1}&\longrightarrow p_i,
  \end{aligned}
\]
where $1 \leqslant i \leqslant n+1,$ 
\[ S_0=S_{n+1}=\{\infty\},\;p_0=p_{n+1}=\infty, \]
and
\[ R_0=K_0=K_1=K_{n+1}=R_{n+1}=\varnothing. \]
Then we can conclude that
\[\sum_{i=1}^n (-1)^i\# S_i \! =\! \sum_{i=1}^n (-1)^i(|K_i|+|R_i|) \! = \! \sum_{i=1}^n (-1)^i(|K_i|+|K_{i+1}|)\!=\!-| K_1 |+(-1)^n |K_{n+1} | \! =\!0. \]
\pre
\cor
\label{equality of e.s. corollary}
Let $N$ be a finitely generated, semipositive binoid and let $I$ be an ideal of $N$. If $J$ is an $N_+$-primary ideal of $N$ then 
\[\# N/J+\# I\cap J/(I+J)=\# I/(I+J)+\# (N/I)/(J+N/I).\]
\core
\pr
We want to apply Proposition \ref{general equation of exact seq} to the strongly exact sequence
\[\infty\longrightarrow I\cap J/(I+J) \longrightarrow I/(I+J) \longrightarrow N/J\longrightarrow (N/I)/(J+N/I)\longrightarrow\infty \]
from Proposition \ref{strongly exact sequence}. To do this we have to show that the involved $N$-sets are finite.
We know that $N/J$ is a finite set.
Also we know that $I$ is a finitely generated $N$-set, so by Proposition \ref{surj map to smash} we have a surjective homomorphism
$(N/J)^{\cupdot r}\rightarrow I\wedge_N N/J$. Hence $|I\wedge_N N/J |=|I/(J+I)|\leqslant |N/J|^r$,
which is a finite set. So we can apply Proposition \ref{general equation of exact seq} and get the result.
\pre

\section{Algebras and Modules}

In this section we will assume that $K$ is a commutative ring (or just a field). We associate to a binoid a binoid algebra over $K$
essentially in the same way how monoids yield monoid algebras.

\df
The \emph{binoid algebra} of a binoid $N$ is the quotient algebra \[ KN/\langle X^\infty\rangle=:K[N], \]
where $KN$ is the monoid algebra of $N$ and $\langle X^\infty\rangle$ is the ideal in $KN$ generated by the element $X^\infty$.
\dfe
So we can consider $K[N]$ as the set of all formal sums $\sum_{m\in M}r_m X^m,$ where $M\subseteq N^\bullet$ is finite, $r_m\in K$ and the multiplication is given by
\[ r_nX^n\cdot r_m X^m=\left\{
    \begin{array}{ccc}
     (r_nr_m)X^{n+m},&\text{ if } n+m\in N^\bullet \\
     0,&\text{ if } n+m=\infty. \\
      \end{array}
       \right. \]
For an $N$-set $(S,p)$ we define the $K[N]$-module $K[S]$ as the set of all formal sums $\sum_{s\in U}r_s X^s,$ where $U\subseteq S^\bullet=S \setminus \{p\}$ is finite, $r_s\in K$
and the multiplication is given by
\[(\sum_{m\in M}r_m X^m)\cdot (\sum_{s\in U}r_s X^s)=\left\{
    \begin{array}{ccc}
     \sum_{m\in M,s\in U}r_mr_s X^{m+s},&\text{ if } m+s\in S^\bullet \\
     0,&\text{ if } m+s=p. \\
      \end{array}
       \right. \]
Here $M$ is a finite subset of $N$ and $U$ is a finite subset of $S$. For an ideal $I \subseteq N$ we get an ideal
\[ K[I]:=\{\sum_{a\in J}r_a X^a\mid  J\subseteq I \text{ finite subset }\} \] of $K[N]$. In this case we have the natural identification
$K[N/I]\cong K[N]/K[I]$. In the same way we have $K[S/T] \cong K[S]/K[T]$. For a finite $N$-set $S$ we have $\# S=\dim_K K[S]$.

\prop
\label{exact sequence of K algebra}
Let $N$ be a binoid, 
\[\infty \longrightarrow S_1 \longrightarrow S_2 \longrightarrow\cdots\longrightarrow S_n\longrightarrow \infty \]
a strongly exact sequence of finite $N$-sets and $K$ a commutative ring. Then we have an exact sequence of $K[N]$-modules 
\[ 0 \longrightarrow K[S_1] \longrightarrow K[S_2] \longrightarrow\cdots\longrightarrow K[S_n]\longrightarrow 0. \]
\prope
\pr
See \cite[Proposition 1.9.6]{Batsukhthesis}.
\pre

\ex
Let $S=\{a,b,p\}$ be as in Example \ref{example of not strongly exact}. Then we have an exact sequence of $N$-sets
\[\infty\longrightarrow\ker f=\infty\xrightarrow{\;i\;\,} S\xrightarrow{\;f\;} S\longrightarrow S/\im f\longrightarrow\infty.\]
We have
\[ K[f](X^a-X^b)=X^{f(a)}-X^{f(b)}=X^b-X^b=0, \]
but $X^a-X^b\notin \im K[i]=\{\infty\}$.
So strong exactness is a necessary condition for Proposition \ref{exact sequence of K algebra}.
\exe

\prop
\label{smash to tensor}
Let $N$ be a binoid and $S,T$ the $N$-sets. Then we have 
\[K[S\wedge_N T]\cong K[S]\otimes_{K[N]} K[T]\] 
and 
\[K[S\cupdot T]\cong K[S]\oplus K[T].\]
\prope
\pr
See Corollary 3.5.2 in \cite{Simone} and \cite[Proposition 1.9.8]{Batsukhthesis}.
\pre

\prop
\label{K algebra of quotient}
Let $N$ be a binoid, $S$ be an $N$-set and $I$ be an ideal of $N$. Then we have \[K[S/(I+S)]\cong K[S]/(K[I]K[S]). \]
\prope
\pr
By Proposition \ref{quotient} we have $S/(I+S)\cong (N/I)\wedge_N S$, so $K[S/(I+S)]\cong K[(N/I)\wedge_N S]$. Hence from the compatibility of the $K$-functor
with the residue class construction, Proposition \ref{smash to tensor} and the general isomorphism $R/I \otimes_R V \cong V/IV$ we get
\[
\begin{aligned}
K[S/(I+S)]&\cong K[(N/I)\wedge_N S] \\
&\cong K[N/I]\otimes_{K[N]} K[S]\\
&\cong (K[N]/K[I])\otimes_{K[N]} K[S] \\
&\cong K[S]/(K[I]K[S]).
\end{aligned} \] 
\qedhere
\pre

\lem
\label{semipositive binoid algebra}
Let $N$ be a semipositive binoid and $K$ some field of characteristic $p$ which does not divide $|N^\times|$. Then $K[N_+]$ is the intersection of finitely many maximal ideals. In particular, if $N$ is a positive binoid then $K[N_+]$ is a maximal ideal of $K[N]$.
\leme
\pr
Because of the fundamental theorem for finite abelian groups we can write
\[ N^\times=\z/(\alpha_1)\times\cdots\times \z/(\alpha_r).\]
By assumption $p$ does not divide $\alpha_1,\dots,\alpha_r$. So from these conditions we can deduce that
\[\begin{aligned}
\overline{K}N^\times &\cong \overline{K}[X_1,\dots,X_r]/(X_1^{\alpha_1}-1,\dots,X_r^{\alpha_r}-1)\\
&=\overline{K}[X_1,\dots,X_r]/\big((X_1-\xi_{11})\cdots(X_1-\xi_{1\alpha_1}),\dots,
(X_r-\xi_{r1})\cdots(X_r-\xi_{r\alpha_r})\big)\\
&\cong\overline{K}[X_1]/\big((X_1-\xi_{11})\cdots(X_1-\xi_{1\alpha_1})\big)\otimes\cdots
\otimes\overline{K}[X_r]/\big((X_r-\xi_{r1})\cdots(X_r-\xi_{r\alpha_r})\big)\\
&\cong\overline{K}^{|N^\times|},
  \end{aligned} \]
where $\overline{K}$ is the algebraic closure of $K$ and $\xi_{ij}$ are the $\alpha_i$-th roots of unity.
Hence the maximal ideals of $\overline{K}N^\times$ have the form 
\[\mathfrak{m}_{i_1,\dots,i_r}=(X_1-\xi_{1i_1},\dots,X_r-\xi_{ri_r}).\]
So we have finitely many maximal ideals in $\overline{K} N^\times$ with this form. 
We also know that the intersection of all maximal ideals of $\overline{K} N^\times$ is equal to $\nil (\overline{K} N^\times)$ and this is $0$. 
Under the homomorphism $KN^\times\hookrightarrow\overline{K}N^\times$ the preimage of a maximal ideal is maximal and therefore the intersection of the maximal ideals of
$KN^\times$ is 0 as well.

Let $K[\pi]:K[N]\rightarrow K[N/N_+]$ be the homomorphism induced by $\pi:N\rightarrow N/N_+\cong (N^\times)^\infty$. Then $K[\pi]^{-1}(\mathfrak{m_i})$ is a maximal ideal
of $K[N]$, where $\mathfrak{m_i}$ is a maximal ideal of $K[N^\times]$. So
\[K[N_+]=K[\pi]^{-1}(0)=\bigcap_i K[\pi]^{-1}(\mathfrak{m_i}).\qedhere \]
\pre

\section{The Hilbert-Kunz function of a binoid}

In this Section we introduce the Hilbert-Kunz function of a binoid. It is defined for a natural number $q$, an $N_+$-primary ideal $\mathfrak{n}$ and a finitely generated $N$-set $T$,
where $N$ fulfills some natural properties. The first lemma ensures that this function exists.

\lem
Let $N$ be a finitely generated, semipositive binoid, $T$ a finitely generated $N$-set and $\mathfrak{n}$ an $N_+$-primary ideal of $N$. Then for every positive integer $q$ we have 
\[\# T/([q]\mathfrak{n}+T)\leqslant r|N^\times|+D q^s,\] where $r$ is the number of generators of $T$, $s$ is the number of generators of $N_+$ and $D$ is some constant.
\leme
\pr
Let $t_1,\dots,t_r$ be generators of $T$ and $n_1,\dots,n_s$ be generators of $N_+$. If $t\in T$ then either \[ t=u+t_i,\] where $u\in N^\times$, or \[t=a_1n_1+\cdots+a_sn_s+t_i,\] where $a_j\in \n,1\leqslant i \leqslant r$.
There are at most $r|N^\times|$ elements of the first type. In the second case we know by the primary property that there exist $d_i\in\n$ such that
$d_in_i\in \mathfrak{n}$ for $1\leqslant i \leqslant s$. So if $a_j\geqslant qd_j$, for some $j$, then $t\in [q]\mathfrak{n}+T$, which means
\[ \# T/([q]\mathfrak{n}+T)\leqslant r|N^\times|+r\cdot q^s\cdot\prod_{i=1}^s d_i.\qedhere \]
\pre

\df
Let $N$ be a finitely generated, semipositive binoid, $T$ a finitely generated $N$-set and $\mathfrak{n}$ an $N_+$-primary ideal of $N$. Then we call the number \[\hkf^N(\mathfrak{n},T,q)=\hkf(\mathfrak{n},T,q):=\# T/([q]\mathfrak{n}+T)=|T/([q]\mathfrak{n}+T)|-1\] the \emph{Hilbert-Kunz function} of $\mathfrak{n}$ on the $N$-set $T$ at $q$.
\dfe
In particular, for $T=N$ and $\mathfrak{n}=N_+$ we have \[\hkf(N,q):=\hkf(N_+,N,q)=\#N/[q]N_+.\]
\ex
\label{N^n}
Let $N=(\n^n)^\infty$ then $\hkf(N,q)=\# N/[q]N_+=q^n$. Since $N_+=(\n^n)^\infty\setminus \{0\}$ we have $[q]N_+=\bigcup_{i=1}^n qe_i+N$, where $e_i$ is the $i$-th standard vector of $\n^n$. Hence \[N/[q]N_+=\{(a_1,\dots,a_n)\in \n^n\mid 0\leqslant a_i\leqslant q-1\}\cup\{\infty\},\] which means $\# N/[q]N_+=q^n$.
\exe

\lem
\label{HKF of annihilator}
Let $N$ be a finitely generated and semipositive binoid, $(T,p)$ a finitely generated $N$-set and $\mathfrak{n}$ an $N_+$-primary ideal of $N$. Let $\mathfrak{a}\subseteq \ann T$ be an ideal of $N$, where $\ann T$ is the annihilator of $T$. 
Then $T$ is also an $N/\mathfrak{a}$-set and we have 
\[\hkf^N(\mathfrak{n},T,q)=\hkf^{N/\mathfrak{a}}((\mathfrak{n}\cup\mathfrak{a})/\mathfrak{a},T,q).\]
\leme
\pr
Let us first show that $(\mathfrak{n}\cup\mathfrak{a})/\mathfrak{a}$ is an $(N/\mathfrak{a})_+$-primary ideal of $N/\mathfrak{a}$.
If there is an element $[m]\in (N/\mathfrak{a})_+$,
where $m\in N$ then there exists $l\in\n$ such that $lm\in\mathfrak{n}$ and so $l[m]=[lm]\in(\mathfrak{n}\cup\mathfrak{a})/\mathfrak{a}$.
Let us define the action $(N/\mathfrak{a})\times T\longrightarrow T$ by $(\overline{n},t)\longrightarrow \overline{n}+t=n+t$.
Then it is easy to see that this is well defined, that 
\[ \begin{tikzpicture}[node distance=2cm, auto]
  \node (N) {$N\times T$};
  \node (T) [right of=N]  {$T$};
  \node (N1) [below of=N]   {$(N/\mathfrak{a})\times T$};
  \node (T1) [below of=T]   {$T$};
  \draw[->] (N) to node {}(T);
  \draw[->] (N) to node {}(N1);
  \draw[->] (N1) to node {}(T1);
  \draw[->] (T) to node {}(T1);
  \end{tikzpicture}
\]
commutes and so $T$ is an $N/\mathfrak{a}$-set. For every $q\in \n$ we have $[q]\mathfrak{n}+T=[q](\mathfrak{n}\cup\mathfrak{a})/\mathfrak{a}+T$,
as this is the image of $[q]\mathfrak{n}\times T$. Therefore
\[ \#T/([q]\mathfrak{n}+T)=\#T/([q](\mathfrak{n}\cup\mathfrak{a})/\mathfrak{a}+T).\qedhere \]
\pre

\lem
\label{HKF ineq}
Let $N$ be a finitely generated, semipositive binoid, $S$ and $T$ finitely generated $N$-sets and $\mathfrak{n}$ an $N_+$-primary ideal of $N$. 
Suppose that we have a surjective $N$-set homomorphism $\phi:S\rightarrow T$. 
Then for all $q$ 
\[\hkf(\mathfrak{n},S,q)\geqslant \hkf(\mathfrak{n},T,q).\]
\leme
\pr
By definition $[q]\mathfrak{n}$ is an ideal of $N$. So by Lemma \ref{surj hom} we have a surjective homomorphism
$S/([q]\mathfrak{n}+S)\longrightarrow T/([q]\mathfrak{n}+T)$. Hence
\[\hkf(\mathfrak{n},S,q)= |S/([q]\mathfrak{n}+S) |-1\geqslant | T/([q]\mathfrak{n}+T) |-1=\hkf(\mathfrak{n},T,q).\qedhere \]
\pre
\df
Let $N$ be a finitely generated, semipositive binoid, $T$ a finitely generated $N$-set and $\mathfrak{n}$ an $N_+$-primary ideal of $N$. 
Then the \emph{Hilbert-Kunz multiplicity} of $\mathfrak{n}$ on the $N$-set $T$ is defined by 
\[ e_{HK}(\mathfrak{n},T):=\lim_{q\rightarrow \infty} \frac{\hkf^N(\mathfrak{n},T,q)}{q^{\dim N}},\] 
if this limit exists.

In particular for $T=N$ and $\mathfrak{n}=N_+$ we set $e_{HK}(N):=e_{HK}(N_+,N)$ and denote it the \emph{Hilbert-Kunz multiplicity} of $M$.
\dfe
Note that here we work with the combinatorial dimension of $N$.

\tr
\label{bound of HKF}
Let $N$ be a finitely generated, semipositive binoid and $\mathfrak{n}$ an $N_+$-primary ideal of $N$. If $T$ is a finitely generated $N$-set and $\hkf(\mathfrak{n},N,q)\leqslant Cq^{\dim N}$ for every $q\in \n$ and some constant $C$ (in particular, if $e_{HK}(\mathfrak{n},N)$ exists) then there exists $\alpha$ such that $\hkf(\mathfrak{n},T,q)\leqslant \alpha q^{\dim N}$.
\tre
\pr
By Lemma \ref{canonical isomorphism}, we have a canonical isomorphism $N^{\cupdot r}/([q]\mathfrak{n}+N^{\cupdot r})\longrightarrow (N/[q]\mathfrak{n})^{\cupdot r}$, which means that $\# ( N^{\cupdot r}/([q]\mathfrak{n}+N^{\cupdot r})) = \# ( (N/[q]\mathfrak{n})^{\cupdot r} ) = r\cdot \#(N/[q]\mathfrak{n})\leqslant r\cdot Cq^{\dim N}$.

From Lemma \ref{surj map for f.g. set}, we have a surjective map $\phi:N^{\cupdot r}\longrightarrow T$, where $r$ is the number of generators of $T$. Also by Lemma \ref{surj hom}, we know that there exists a surjective homomorphism 
\[\tilde{\phi}:N^{\cupdot r}/([q]\mathfrak{n}+N^{\cupdot r}) \longrightarrow T/([q]\mathfrak{n}+T).\]  
So by Lemma \ref{HKF ineq} we have $\hkf(\mathfrak{n},T,q)\leqslant \alpha q^{\dim N}$, where $\alpha=r\cdot C$.
\pre

\lem
\label{multiplicity of smash product}
Let $M$ and $N$ be finitely generated, semipositive binoids. Then we have 
\[\hkf(M\wedge N,q)=\hkf(M,q)\cdot \hkf(N,q).\]
\leme
\pr
We have
\[(M\wedge N)_+=(M\wedge N)\setminus (M^\times\wedge N^\times)=(M_+\wedge N)\cup (M\wedge N_+). \]
Take an element $q(m\wedge n)+m_1\wedge n_1\in [q](M\wedge N)_+$, where $m\wedge n\in (M\wedge N)_+$. Then $m\wedge n\in M_+\wedge N$ or $m\wedge n\in M\wedge N_+$, so $qm+m_1\in [q]M_+$ or $qn+n_1\in [q]N_+$, which means 
\[[q](M\wedge N)_+\subseteq ([q]M_+\wedge N)\cup (M\wedge [q]N_+).\]
We can also easily check the other inclusion, so we have 
\[[q](M\wedge N)_+= ([q]M_+\wedge N)\cup (M\wedge [q]N_+).\]
Hence, from this result and Lemma \ref{smash of quotients}, we get 
\[(M\wedge N)/[q](M\wedge N)_+\cong (M/[q]M_+)\wedge (N/[q]N_+).\] 
By assumption $M,N$
are semipositive, so we know $M/[q]M_+$ and $N/[q]N_+$ are finite binoids. Hence, by Proposition \ref{number of elements of smash}, we can conclude that
\[ \# \big((M\wedge N)/[q](M\wedge N)_+\big)=\# (M/[q]M_+) \cdot \# (N/[q]N_+).\qedhere \]
\pre

\tr
\label{multiplicity of smash}
Let $M$ and $N$ be binoids such that $e_{HK}(M)$ and $e_{HK}(N)$ exist. Then $e_{HK}(M\wedge N)$ exists and 
\[e_{HK}(M\wedge N)=e_{HK}(M)\cdot e_{HK}(N).\]
\tre
\pr
By definition
\begin{align*}
e_{HK}(M)\cdot e_{HK}(N)&=\lim_{q\rightarrow \infty} \frac{\hkf(M,q)}{q^{\dim M}}\cdot \lim_{q\rightarrow \infty} \frac{\hkf(N,q)}{q^{\dim N}}\\
&=\lim_{q\rightarrow \infty} \frac{\hkf(M,q)\cdot \hkf(N,q)}{q^{\dim M+\dim N}}\\
&\overset{\text{Lemma } \ref{multiplicity of smash product}}{=}\lim_{q\rightarrow \infty} \frac{\hkf(M\wedge N,q)}{q^{\dim M+\dim N}}\\
&\overset{\text{Theorem } \ref{dimension of smash product}}{=}\lim_{q\rightarrow \infty} \frac{\hkf(M\wedge N,q)}{q^{\dim M\wedge N}}\\
&=e_{HK}(M\wedge N).
\end{align*}
\pre

\section{Reduction to the integral and reduced case}

We want to show that the Hilbert-Kunz multiplicity of a finitely generated semipositive cancellative reduced binoid exists by reducing it to the toric (not necessarily normal) case, i.e. to the case of a integral cancellative torsionfree binoid which was proven by Eto (\cite{Eto}). If $N$ is such a monoid and $N \subseteq \hat{N} \subseteq \diff N$ its normalization,
and $   {\mathfrak n} $ a primary ideal with generators $f_1, \ldots , f_n$,
then the Hilbert-Kunz multiplicity $e_{HK} ( {\mathfrak n},N)$ equals $ e_{HK} (  {\mathfrak n} \hat{N}, \hat{N} )$, where $ {\mathfrak n} \hat{N}$ denotes
the extended ideal (see Lemma \ref{finite birational extension} below). If we write the normal toric positive monoid $\hat{N}$ as $\hat{N} =C \cap \Gamma $ with a polyhedral cone $C$ and the lattice $\Gamma=\diff N \setminus \{\infty\}$,
then the Hilbert-Kunz multiplicity is the volume of the region inside $C$ and outside of $ \bigcup_{i =1}^n f_i +C$.

In this Section we show how to reduce the existence and rationality of the Hilbert-Kunz multiplicity to the integral and the reduced case.
\lem
\label{reduction of torsion free}
Suppose that for all finitely generated, semipositive, cancellative, integral and torsion-free binoids of dimension less or equal to $d$ the Hilbert-Kunz multiplicity exists. Then
\begin{enumerate}
\item Let $N$ be a finitely generated, semipositive, cancellative, torsion-free, reduced  binoid of dimension $d$ and let $\{\mathfrak{p}_1,\dots,\mathfrak{p}_s\}$ be all minimal prime ideals of  dimension $d$. Then the Hilbert-Kunz multiplicity of $N$ exists, and 
\[e_{HK}(N)=\sum_{i=1}^s e_{HK}(N/\mathfrak{p}_i).\]
\item Let $N$ be a finitely generated, semipositive, cancellative binoid of dimension $d$ and torsion-free except for nilpotent elements. Then we have 
\[\hkf(N,q)\leqslant Cq^{\dim N},\] 
for all $q$, where $C$ is some constant.
\end{enumerate}
\leme
\pr
We have finitely many minimal prime ideals $\{\mathfrak{p}_1,\dots,\mathfrak{p}_n\}$ and $n\geqslant s$.
From reducedness we have $\bigcap_{i=1}^n \mathfrak{p}_i=\nil(N)=\{\infty\}$ by  \cite[Corollary 2.3.10]{Simone}. If $f\notin [q]N_+$ and $f\in [q]N_+\cup \mathfrak{p}_i$ for all $i$,
then $f\in \bigcap_{i=1}^n \mathfrak{p}_i=\{\infty\}$, which implies $N\setminus[q]N_+\subseteq \bigcup_{i=1}^n N\setminus([q]N_+\cup\mathfrak{p}_i)$.
Also by Proposition \ref{quotient to union} and the set identifications
\[ N\setminus[q]N_+\cong (N/[q]N_+)\setminus \{\infty\}, \]
\[ N\setminus([q]N_+\cup\mathfrak{p}_i)\cong (N/([q]N_+\cup\mathfrak{p}_i))\setminus \{\infty\}\]
we can conclude that
\[N\setminus[q]N_+\subseteq \bigcup_{i=1}^n (N/\mathfrak{p}_i)\setminus[q](N/\mathfrak{p}_i)_+. \]
Since $N/\mathfrak{p}_i$ is finitely generated, semipositive, cancellative, torsion-free and integral by assumption we know that 
\[\hkf(N/\mathfrak{p}_i,q)=e_{HK}(N/\mathfrak{p}_i)q^{d_i}+O(q^{d_i-1}),\]
where $d_i:=\dim N/\mathfrak{p}_i$. Hence only the minimal prime ideals 
with dimension $d_i=d$ are important and we can write 
\begin{equation}
\# N/[q]N_+\leqslant \sum_{i=1}^s e_{HK}(N/\mathfrak{p}_i)q^d+O(q^{d-1}).
\end{equation}

Let us denote $H:=N/[q]N_+$. Then 
\[N/([q]N_+\cup\mathfrak{p}_i)=(H\setminus \mathfrak{p}_i)\cup\{\infty\}\] 
and
\[ N/([q]N_+\cup\mathfrak{p}_i\cup\mathfrak{p}_j)=\big((H\setminus \mathfrak{p}_i)\cap (H\setminus \mathfrak{p}_j)\big)\cup\{\infty\}.\]
Hence from set theory we know that
\[ |H\setminus \{\infty\}|\geqslant |\bigcup_{i=1}^n H\setminus\mathfrak{p}_i|\geqslant \sum_{i=1}^n |H\setminus\mathfrak{p}_i|-\sum_{1\leqslant i<j\leqslant n} |(H\setminus \mathfrak{p}_i)\cap (H\setminus \mathfrak{p}_j)|.\]
So from here we have
\[\#  N/[q]N_+\geqslant \sum_{i=1}^n \# N/([q]N_+\cup\mathfrak{p}_i)-\sum_{1\leqslant i<j\leqslant n}\# N/([q]N_+\cup\mathfrak{p}_i\cup\mathfrak{p}_j). \]
But we know that $N/(\mathfrak{p}_i\cup\mathfrak{p}_j)$ is finitely generated, semipositive, cancellative, torsion-free and integral so by Proposition \ref{quotient to union} and the assumption we get 
\begin{align*}
\sum_{i\neq j}\# N/([q]N_+\cup\mathfrak{p}_i\cup\mathfrak{p}_j)&=\sum_{i\neq j}\hkf(N/(\mathfrak{p}_i\cup\mathfrak{p}_j),q)\\
&=\sum_{i\neq j}\Big(e_{HK}(N/(\mathfrak{p}_i\cup\mathfrak{p}_j))q^{d_{ij}}+O(q^{d_{ij}-1})\Big), 
\end{align*}
 where $d_{ij}$ is the dimension of $N/(\mathfrak{p}_i\cup\mathfrak{p}_j)$. By Lemma \ref{dimension properties} we also know that $d_{ij}<\min\{d_i,d_j\}$. So we have 
\begin{equation}
\#  N/[q]N_+\geqslant \sum_{i=1}^s e_{HK}(N/\mathfrak{p}_i)q^d+O(q^{d-1}). 
\end{equation}
Hence from (1) and (2) we get
\[e_{HK}(N)=\lim_{q\rightarrow \infty} \dfrac{\hkf(N,q)}{q^d}=\sum_{i=1}^s e_{HK}(N/\mathfrak{p}_i). \]
To prove the second statement, note that $\nil(N)$ is an ideal of $N$ and $N_{\red}:=N/\nil(N)$ is reduced.
Let $a_1,\dots,a_s$ be the generators of $\nil(N)$ and $k_ia_i=\infty$ but $(k_i-1)a_i\neq \infty$, where $k_i\in \n,1\leqslant i\leqslant s$.
We have a finite decreasing sequence
\[N=M_0\longrightarrow M_1\longrightarrow \cdots \longrightarrow M_{t-1}\longrightarrow M_t=N_{\red},\]
where $M_{k_0+\cdots+k_i+j}=M_{k_0+\cdots+k_i+j-1}/((k_{i+1}-j)a_{i+1}),1\leqslant j\leqslant k_{i+1}-1,0\leqslant i\leqslant s-1$ and $k_0=0$.
Here we have $2(k_{i+1}-j)a_{i+1}=\infty$ in $M_{k_0+\cdots+k_i+j}$.
So in particular there exists a sequence such that 
\[M_{i+1}=M_i/(f_i),2f_i=\infty \text{ in }M_i,0\leqslant i\leqslant t-1\]
and $M_0=N,M_t=N_{\red}$.
Hence there exists also such a sequence with this property of minimal length $l$.
We will use induction on $l$ to prove this Lemma.

For $l=0$ we are in the reduced situation and the statement follows from the case above. So suppose that $l$ is arbitrary and that for smaller $l$ the statement is already proven. Suppose
that a sequence as described is given. This means in particular that we have a strongly exact sequence 
\[\infty \!  \longrightarrow (f)+N  \longhookrightarrow N  \longtwoheadrightarrow N/(f) \! \longrightarrow \! \infty,\]
where $f=f_0$.  Hence $M:=N/(f)$ has a decreasing sequence with the described property of length smaller than $l$, 
so by the induction hypothesis we know that $\# M/[q]M_+\leqslant Dq^d$ for some constant $D$. By Corollary \ref{equality of e.s. corollary} we also have 
\[\# N/[q]N_+ +\# \big((f)\cap [q]N_+\big)/((f)+[q]N_+)=\# ((f)+N)/((f)+[q]N_+)+\# M/([q]N_++M).\]
Since $2f=\infty$, we know that $f$ annihilates the $N$-set $(f)+N$. So $(f)+N$ is a finitely generated $M$-set ($f$ is the only $M$-set generator of $(f)+N).$
Hence by Lemma \ref{surj map for f.g. set} we have a surjective homomorphism $M\longrightarrow (f)+N$, which implies 
\[\# ((f)+N)/((f)+[q]N_+)\leqslant\# M/[q]M_+.\] 
But we also know that $[q]M_+\subseteq [q]N_++M$, so we have that 
\[\# M/([q]N_++M)\leqslant\# M/[q]M_+.\] 
From here we can conclude
that
\[\# N/[q]N_+\leqslant\# M/[q]M_++\# M/([q]N_++M)\leqslant 2Dq^d.\qedhere \]
\pre

The same reduction steps hold when we allow torsion.
\lem
\label{reduction}
Suppose that for all finitely generated, semipositive, cancellative and integral binoids of dimension less or equal to $d$ the Hilbert-Kunz multiplicity exists. Then
\begin{enumerate}
 \item Let $N$ be a finitely generated, semipositive, cancellative and reduced  binoid of dimension $d$ and let $\{\mathfrak{p}_1,\dots,\mathfrak{p}_s\}$ be all minimal prime ideals of dimension $d$. 
 Then the Hilbert-Kunz multiplicity of $N$ exists, and 
\[e_{HK}(N)=\sum_{i=1}^s e_{HK}(N/\mathfrak{p}_i).\]
 \item Let $N$ be a finitely generated, cancellative and semipositive binoid of dimension $d$. Then 
 \[\hkf(N,q)\leqslant Cq^d,\] 
 for all $q$, where $C$ is some constant.
\end{enumerate}
\leme
\pr
The proof is similar to the one of Lemma \ref{reduction of torsion free}.
\pre

The following Lemma reduces in particular the case of a non-normal toric binoid to the normal toric case. A finite $N$-binoid $M$ is an $N$-binoid which is finitely generated as an $N$-set.
\lem
\label{finite birational extension}
Suppose that for all finitely generated, cancellative, semipositive, integral binoids of dimension less than $d$, the Hilbert-Kunz function is bounded by $Cq^{d-1}$ for some constant $C$. Let $N$ be a finitely generated, semipositive, cancellative, integral binoid of dimension $d$ and let $\mathfrak{n}$ be an $N_+$-primary ideal of $N$. Let $M$ be a finite $N$-binoid which is birational over $N$. Suppose that $e_{HK}(\mathfrak{n}+M,M)$ exists. Then $e_{HK}(\mathfrak{n},N)$ exists and it is equal to $e_{HK}(\mathfrak{n}+M,M)$.
\leme
\pr
First note that $\mathfrak{n}+M$ is an $M_+$-primary ideal. By Proposition \ref{strongly exact sequence} applied to the $N$-sets $S=N,T=M$ and the ideal $J=[q]\mathfrak{n}$ we have the strongly exact sequence 
\[\infty \! \longrightarrow (N\cap([q]\mathfrak{n}+M))/[q]\mathfrak{n} \! \longrightarrow \! N/[q]\mathfrak{n} \! \longrightarrow M/([q]\mathfrak{n}+M) \! \longrightarrow  (M/N)/([q]\mathfrak{n}+M/N) \! \longrightarrow\infty.\]
By Proposition \ref{general equation of exact seq},
we know that
\[ \# (N\cap([q]\mathfrak{n}+M))/[q]\mathfrak{n} + \# M/([q]\mathfrak{n}+M) = \# N/[q]\mathfrak{n} +\#(M/N)/([q]\mathfrak{n}+M/N) \]
and we have 
\[ \# M/([q]\mathfrak{n}+M) \leqslant \# N/[q]\mathfrak{n} +\#(M/N)/([q]\mathfrak{n}+M/N). \]
It is easy to check that there exists (a common denominator) $b\in N$ such that $b+M\subseteq N$, and that $I:=b+M\neq\{\infty\}$ is an ideal of $N$ which is isomorphic to $M$, because $N$ is an integral binoid. We also know that $I$ annihilates $M/N$,
so $M/N$ is an $N/I$-set and $N/I$ is a finitely generated, semipositive, cancellative binoid
and torsion-free up to nilpotence.
So by Lemma \ref{HKF of annihilator} we have that
\[\# (M/N)/([q]\mathfrak{n}+M/N)=\hkf(\mathfrak{n},M/N,q)=\hkf^{N/I}((I\cup\mathfrak{n})/I,M/N,q)\]
and by Lemma \ref{dimension properties} we know that $\dim N/I<d$, so from the assumptions and by Lemma \ref{reduction of torsion free} we get 
\[\hkf^{N/I}((I\cup\mathfrak{n})/I,M/N,q)\leqslant C q^{d-1}.\]
Hence we have
\begin{equation}
 \#  M/([q]\mathfrak{n}+M)-Cq^{d-1} \leqslant \# N/[q]\mathfrak{n}.
\end{equation}
By Corollary \ref{equality of e.s. corollary}, we deduce
\[\# N/[q]\mathfrak{n}\leqslant\# I/(I+[q]\mathfrak{n})+\# (N/I)/([q]\mathfrak{n}+N/I).\]
But we know that $\# (N/I)/([q]\mathfrak{n}+N/I)=\hkf(\mathfrak{n},N/I,q)$ and $I$ annihilates $N/I$. Similarly we can show that 
\[\hkf^N(\mathfrak{n},N/I,q)=\hkf^{N/I}((I\cup\mathfrak{n})/I,N/I,q)\leqslant C'q^{d-1},\] 
where $C'$ is some constant.
Since $I\cong M$, we have that $\# I/(I+[q]\mathfrak{n})=\# M/([q]\mathfrak{n}+M)$, so we can conclude that 
\[\# N/[q]\mathfrak{n}\leqslant\# M/([q]\mathfrak{n}+M)+C'q^{d-1}.\]
From the last result and $(3)$ we have that
\[e_{HK}(\mathfrak{n},N)=  \operatorname{lim}_{q \rightarrow \infty}      \frac{  \# M/([q]\mathfrak{n}+M)}{q^d}.\]
Since extended ideals commute with Frobenius sums, we have $\# M/([q]\mathfrak{n}+M)=\# M/[q](\mathfrak{n}+M)$, which means that 
\[e_{HK}(\mathfrak{n},N)=e_{HK}(\mathfrak{n}+M,M).\]
\pre

\tr
\label{not integral, reduced}
Let $N$ be a finitely generated, positive, cancellative, torsion-free and  reduced binoid and let $\{\mathfrak{p}_1,\dots,\mathfrak{p}_s\}$ be all minimal prime ideals with dimension $d=\dim N$.
Then 
\[e_{HK}(N)=\lim_{q\rightarrow \infty} \dfrac{\hkf(N,q)}{q^d}=\sum_{i=1}^s e_{HK}(N/\mathfrak{p}_i)\] 
exists and it is a rational number.
\tre
\pr
This follows from Lemma \ref{reduction of torsion free} (1) and the toric case.
\pre

\ex
\label{simplicial}
Let $\trian$ be a simplicial complex on the vertex set $V$. The binoid associated to $\trian$ is given by $F(V)/I_{\trian}=:M_{\trian}$, where $I_{\trian}$ is
the ideal $\{ f\in F(V)\mid  \supp(f)\nsubseteq \trian\}$ of the free binoid $F(V)\cong (\n^{|V|})^\infty$. In this case we have
$e_{HK}(M_{\trian})=k $,
where $k$ is the number of facets of maximal dimension. This rests on Theorem \ref{not integral, reduced} and the facts that simplicial binoids are positive cancellative torsion-free and reduced,
that faces correspond to prime ideals and facets to minimal prime ideals and that $M_{\trian}/ \mathfrak{p} \cong {\mathbb N}^{\operatorname{dim} M_{\trian} }$
for the minimal prime ideals, see \cite[Chapter 6]{Simone}.
\exe

\lem
\label{not integral and torsion free up to nilpotence}
Let $N$ be a finitely generated, positive, cancellative binoid which is torsion-free up to nilpotence. Then $\hkf(N,q)\leqslant Cq^d$,where $C$ is some constant and $d=\dim N$.
\leme
\pr
This follows from the toric case and Lemma \ref{reduction of torsion free} (2).
\pre
\cor
\label{bound for N/I}
Let $N$ be a finitely generated, integral, positive, cancellative and torsion-free binoid and let $I\neq \infty$ be an ideal of $N$. Then $\hkf(N/I,q)\leqslant Cq^d$,where $C$ is some constant and $d=\dim N/I$.
\core
\pr
We know by \cite[Lemma 2.1.20]{Simone}, that $N/I$ is torsion-free up to nilpotence and finitely generated, positive, cancellative.
So by Lemma \ref{not integral and torsion free up to nilpotence} we have the result.
\pre

\lem
\label{ideal case}
Let $N$ be a finitely generated, positive, cancellative, integral and torsion-free binoid and let $I\neq \infty$
be an ideal of $N$. Then $e_{HK}(I)$ exists and is equal to $e_{HK}(N)$.
\leme
\pr
By Corollary \ref{equality of e.s. corollary} we have
\[  \# N/[q]N_+ +\# I\cap [q]N_+/(I+[q]N_+)=\# I/(I+[q]N_+)+\# (N/I)/[q](N/I)_+ . \]
By Corollary \ref{bound for N/I} we know that $\hkf(N/I,q)\leqslant Dq^{d'}$, where $d'=\dim N/I<d$. Hence from $(2.11)$ we get
\begin{equation}
\# I/(I+[q]N_+)\geqslant \hkf(N,q)-\# (N/I)/[q](N/I)_+\geqslant \hkf(N,q)-Dq^{d'}.
\end{equation}

Let $\infty\neq f\in I$. Then we have a strongly exact sequence $\infty\rightarrow (f)\rightarrow I\rightarrow I/(f)\rightarrow\infty$. So we have 
\[\# I/(I+[q]N_+)\leqslant \# (f)/((f)+[q]N_+)+\# (I/(f))/[q](I/(f))_+.\] 
By Lemma \ref{HKF of annihilator}, the $N$-set  $I/(f)$ is an $N/(f)$-set.
By Corollary \ref{bound for N/I} we know that $\hkf(N/(f),q)\leqslant Dq^{d''}$, where $d''=\dim N/(f)<d$ and by Theorem \ref{bound of HKF} we have $\hkf(I/(f),q)\leqslant \alpha q^{d''}$. Hence we can conclude that
\begin{equation}
 \# I/(I+[q]N_+)\leqslant \# (f)/((f)+[q]N_+)+\alpha q^{d''}.
\end{equation}

But we know that $(f)=f+N$ is isomorphic to $N$ as an $N$-set so 
\[\# (f)/((f)+[q]N_+)=\hkf(N,q).\]
Now by $(4)$ and $(5)$ we have 
\[\begin{aligned}e_{HK}(N)&=\lim_{q\rightarrow\infty} \frac{\hkf(N,q)-Dq^{d'}}{q^d}\\
&\leqslant \lim_{q\rightarrow\infty} \frac{\# I/(I+[q]N_+)}{q^d}=e_{HK}(I)\\
&\leqslant\lim_{q\rightarrow\infty} \frac{\hkf(N,q)+\alpha q^{d''}}{q^d}=e_{HK}(N)\end{aligned}\]
\pre

\section{Integral and cancellative binoids with torsion}

Let $N$ be a finitely generated, semipositive, cancellative and integral binoid. We know that
\[N\subseteq \diff N\cong(\z^m\times T)^\infty=(\z^m)^\infty\wedge T^\infty, \]
where smashing is over the trivial binoid $\t=\{0,\infty\}$.
Here  $T$ is the torsion part of the difference group, which is a finite commutative group, hence $T=\z/k_1\times\cdots \times \z/k_l$. We will write elements $x\in N^\bullet$ as $x=f\wedge t$ with $f\in\z^m$ and $t\in T$. This representation is unique.

The relation $\sim_{\tf}$ on $N$ given by $a\sim_{\tf} b$ if $na=nb$ for some $n\geqslant 1$ is a congruence and $N_{\tf}:=N/\sim_{\tf}$ is a torsion-free binoid which we call the \emph{torsion-freefication} of $N$. If $F=\{f\in (\z^m)^\infty\mid  \exists t\in T, f\wedge t\in N\}$ then $N_{\tf}\cong F\subseteq \z^m$. Hence we may assume $N\subseteq F\wedge T^\infty$ and $f\wedge t_1\sim_{\tf} g\wedge t_2$ if and only if there exists $n\in\n$ such that $n(f\wedge t_1)=n (g\wedge t_2)$, which means $f=g\in F$. We define the subsets $F\wedge t:=\{g\wedge t \mid g\in F\}$, $F_t:=\{f\in F \mid f\wedge t\in M\}$ and $M_t:=F_t\wedge t$.

\prop
Let $F$ be a finitely generated, positive, cancellative, integral and torsion-free binoid. If $T$ is a finite group then
$e_{HK}(F\wedge T^\infty)=e_{HK}(F)\cdot |T|$.
\prope
\pr
This follows from the toric case, Theorem \ref{multiplicity of smash} and the fact that for a finite group binoid the Hilbert Kunz multiplicity is just the order of the group.
\pre

In the following we write $N\subseteq F\wedge T^\infty$, where $F\cong N_{\tf}$, $T$ is a finite group and $\diff N=\diff F\wedge T^\infty$. The strategy is to relate the Hilbert-Kunz multiplicity of $N$ with that of $F\wedge T^\infty$.

\lem
\label{correspondence of ideals}
Let $F$ be a binoid and $T^\infty$ be a group binoid. Then we have a bijection between ideals of $F$ and ideals of $F\wedge T^\infty$. The $F_+$-primary ideals correspond to $(F\wedge T^\infty)_+$-primary ideals.
\leme
\pr
We have an inclusion 
\[F \xrightarrow{\;i\;} F\wedge T^\infty,~f \longmapsto f\wedge 0.\] 
So we can consider for an ideal in $F$ its extended ideal in $F\wedge T^\infty$, in other words we use the map $\mathfrak{a}\mapsto \mathfrak{a}+F\wedge T^\infty$. 
For ideal generators $f_j\wedge t_j,j\in J$, in $F\wedge T^\infty$ we have 
\[\langle f_j\wedge t_j\mid j\in J\rangle=\langle f_j\wedge 0\mid j\in J\rangle\subseteq F\wedge T^\infty,\] 
because $0\wedge t_j$ are units.
Hence the map $\mathfrak{b}\mapsto i^{-1}(\mathfrak{b})$, where $\mathfrak{b}$ is an ideal of $F\wedge T^\infty$, is inverse to the extension map.
Let $\mathfrak{p}\subseteq F$ be an $F_+$-primary ideal and $f\wedge t\in (F\wedge T^\infty)_+$. Then there exists $l\in\n$ such that $lf \in \mathfrak{p}$, so $l(f\wedge t)=lf\wedge lt=(lf\wedge 0)+(0\wedge lt) \in \mathfrak{p}+F\wedge T^\infty$, because $0\wedge lt$ is a unit.
Now let $\mathfrak{q}$ be an $F\wedge T^\infty_+$-primary ideal and $f\in F_+$. Then $f\wedge 0\in F\wedge T^\infty$ is not a unit, so there exists $m\in \n$ such that $m(f\wedge 0)=mf\wedge 0\in \mathfrak{q}$. Hence $mf \in i^{-1}(\mathfrak{q})$.
\pre

\lem
\label{finite birational ext}
Let $N$ be a finitely generated, semipositive, cancellative, integral binoid. If $F=N_{\tf}$ and $T$ is the torsion subgroup of $\diff N$ then $F\wedge T^\infty$ is finite and birational over $N$.
\leme
\pr
We can assume that $\diff N=(\z^d)^\infty\wedge T^\infty$, where $d=\dim N$, $T=\{t_1,\dots,t_m\}$ is a finite abelian group. We know $F\wedge T^\infty \subseteq \diff N$, so it is clear that $F\wedge T^\infty$ is birational over $N$. If $f\wedge t\in (F\wedge T^\infty)^\bullet$ then there exists $t'\in T$ such that $f\wedge t'\in N$ and
\begin{equation}
m(f\wedge t)=mf\wedge mt=mf\wedge 0=mf\wedge mt'=m(f\wedge t')\in N,
\end{equation}
which means that these elements satisfy a pure integral equation over $N$. Let $f_i\wedge t_i,1\leqslant i \leqslant k$, be the generators of $N$, then $\{f_i\wedge t_j,\mid 1\leqslant i \leqslant k,1\leqslant j \leqslant m\}$ will give us an $N$-generating system of $F\wedge T^\infty$.
This means that $F\wedge T^\infty$ is finite over $N$.
\pre
Note that the notion birational makes sense though the corresponding binoid algebras are not integral in general.
\lem
\label{F smash T}
Let $F$ be a finitely generated, positive, cancellative, torsion-free, integral binoid, $T^\infty$ be a torsion group binoid and $\mathfrak{p}$ be an $F_+$-primary ideal of $F$. Then 
\[\hkf^{F\wedge T^\infty}(\mathfrak{p}+F\wedge T^\infty,F\wedge T^\infty,q)=|T|\cdot \hkf^F(\mathfrak{p},F,q).\]
In particular, the Hilbert-Kunz multiplicity exists and 
\[e_{HK}^{F\wedge T^\infty}(\mathfrak{p}+F\wedge T^\infty,F\wedge T^\infty)=|T|\cdot e_{HK}^F(\mathfrak{p},F).\]
\leme
\pr
For every $q\in\n$ we have by Lemma \ref{smash of quotients}
\[F\wedge T^\infty/([q]\mathfrak{p}\wedge T^\infty\cup F\wedge (\infty))=F\wedge T^\infty/[q]\mathfrak{p}\wedge T^\infty\cong F/[q]\mathfrak{p}\wedge T^\infty,\]
and $[q]\mathfrak{p}\wedge T^\infty=[q]\mathfrak{p}+F\wedge T^\infty$. Hence
\[ F\wedge T^\infty/([q]\mathfrak{p}+F\wedge T^\infty)\cong F/[q]\mathfrak{p}\wedge T^\infty,\]
which means
\[ \hkf^{F\wedge T^\infty}(\mathfrak{p}+F\wedge T^\infty,F\wedge T^\infty,q)=|T|\cdot \hkf^F(\mathfrak{p},F,q). \]
The existence of the Hilbert-Kunz multiplicity follows from the toric case.
\pre

Let $N$ be a binoid with $N\subseteq F\wedge T^\infty$, where $F\cong N_{\tf}$, $T$ a finite group. Then we have the following diagram.
\[\begin{tikzpicture}[node distance=2cm, auto]
  \node (F) {F};
  \node (FT) [below of=F]  {$F\wedge T^\infty$};
  \node (N) [left of=FT]   {$N$};
  \draw[->] (N) to node {} (FT);
  \draw[->] (F) to node [swap] {$i$} (FT);
\end{tikzpicture}
\]
where $i(f)=f\wedge 0$.

\tr
\label{HKM with torsion binoid}
Let $N$ be a finitely generated, semipositive, cancellative, integral binoid and $\mathfrak{n}$ be an $N_+$-primary ideal of $N$. Then 
\[e_{HK}^N(\mathfrak{n},N)=|T|\cdot e_{HK}^F(\mathfrak{m},F),\] 
where $F=N_{\tf}$, $T$ is the torsion subgroup of $\diff N$ and $\mathfrak{m}=i^{-1}(\mathfrak{n}+F\wedge T^\infty)$ is an ideal of $F$.
\tre
\pr
We know that $\mathfrak{n}+F\wedge T^\infty$ is a primary ideal and by Lemma \ref{correspondence of ideals} that $\mathfrak{m}$ is an $F_+$-primary ideal of $F$. So by Lemma \ref{F smash T} we know that $e_{HK}(\mathfrak{m}+F\wedge T^\infty,F\wedge T^\infty)$ exists and is equal to $|T|\cdot e_{HK}^F(\mathfrak{m},F)$.
By Lemma \ref{correspondence of ideals} we know that $\mathfrak{m}+F\wedge T^\infty=\mathfrak{n}+F\wedge T^\infty$ so
\[ e_{HK}^{F\wedge T^\infty}(\mathfrak{n}+F\wedge T^\infty,F\wedge T^\infty)=e_{HK}^{F\wedge T^\infty}(\mathfrak{m}+F\wedge T^\infty,F\wedge T^\infty). \]
Because of Lemma \ref{finite birational ext} and using induction over the dimension we can apply Lemma \ref{finite birational extension}.
Hence $e_{HK}^N(\mathfrak{n},N)$ exists and is equal to
\[ e_{HK}^{F\wedge T^\infty}(\mathfrak{n}+F\wedge T^\infty,F\wedge T^\infty)=|T|\cdot e_{HK}^F(\mathfrak{m},F). \]
\pre

\tr
\label{f.g,s.p,c,r binoid}
Let $N$ be a finitely generated, semipositive, cancellative, reduced binoid and $\mathfrak{n}$ be an $N_+$-primary ideal of $N$. Then $e_{HK}(\mathfrak{n},N)$ exists and is rational number.
\tre
\pr
This follows directly from Lemma \ref{reduction} and Theorem \ref{HKM with torsion binoid}.
\pre

\ex
The binoid $\langle x,y\rangle/ax=ay$ (for $a\in \n_+$) can be realized as $\langle (1\wedge 0),(1\wedge 1)\rangle\subseteq \n\wedge (\z/a)^\infty$.
In this case, the torsion-freefication is $\n$ and the Hilbert-Kunz multiplicity is $a$ by Theorem \ref{HKM with torsion binoid},
since $e_{HK}((\n)^\infty)=1$ by Example \ref{N^n} and $HKF( (\z/a)^\infty,q)=|(\z/a)|$.
\exe

\ex
Let $a=(2,1),b=(3,0)\in (\n \times \z/2)^\infty$ be the generators of a binoid $N$. This binoid is not torsion-free, since $(6,1) \neq (6,0)$, but $2 (6,1) \neq 2 (6,0)$.
The binoid is positive, its normalization $\hat{N} = (\n \times \z/2)^\infty$ is only semipositive.
It is not difficult to see that  $N_{\tf}\cong\n^\infty\setminus \{1\}$ and the torsion group is $T=\z/2$. So by Theorem \ref{HKM with torsion binoid},
we have 
\[e_{HK}(N)=|T|\cdot e_{HK}(N_{\tf})=2\cdot 2=4.\]
\exe

\ex
Let $N=\langle X,Y,Z\rangle/4X+12Y=16Z$ be a binoid. From \cite[Lemma 2.2.9]{Batsukhthesis} we have an injective binoid homomorphism $\phi:N\rightarrow \big(\z^2\times \z/16\z\big)^\infty$.
So $\phi(N)$ has generators $(16,0,0), (0,16,0), (4,12,1)$. If we choose the new generators $u=(4,-4,1),~v=(0,16,0),~w=(0,0,4)$, then $\phi(X)=4u+v-w$, $\phi(Y)=v$, $\phi(Z)=u+v$.
Hence the difference group of our binoid is isomorphic to $\z^2\times \z/4$.
The torsion-freefication of $N$ is $F=\langle X',Y',Z'\rangle/X'+3Y'=4Z'$ with $e_{HK}(F)=\frac{13}{4}$ by the toric case,
so $e_{HK}(N)=13$ by Theorem \ref{HKM with torsion binoid}.
\exe

\section{Hilbert-Kunz function of binoid rings}

We finally want to relate the Hilbert-Kunz function of a binoid $N$ to the Hilbert-Kunz function of its binoid algebra $K[N]$. However, $K[N]$ is not a local ring. If $N$ is positive, then $K[N]$ contains the unique combinatorial maximal ideal $K[N_+]$ and we work with the localization $K[N]_{K[N_+]}$.
In this setting we get by counting the dimension and Proposition \ref{K algebra of quotient}, immediately
\[ HKF^N\!(\mathfrak{n},S,q) \! =\! \# \big(S/(S+[q]\mathfrak{n})\big) \!=\!\dim_K\! K[S]/(K[S]\cdot K[\mathfrak{n}]^{[q]}) \! =\! HKF^{K[N]}(K[\mathfrak{n}],K[S],q). \]
So the rationality results of the previous sections translates directly to results on the Hilbert-Kunz multiplicity of the localization of a binoid algebra.

This translation is more involved for semipositive binoids. We have seen in Lemma \ref{semipositive binoid algebra}
that $K[N_+]$ is the intersection of finitely many maximal ideals $\mathfrak{m}_1,\dots,\mathfrak{m}_r$.
Then $T=K[N]\setminus\mathfrak{m}_1\cap\cdots\cap K[N]\setminus\mathfrak{m}_r$
is a multiplicatively closed subset of $K[N]$ and the localization $K[N]_T$ is a semilocal ring with maximal ideals
$\mathfrak{m}_i \cdot K[N]_T$ and $K[N_+]\cdot K[N]_T=\bigcap_{i=1}^r \mathfrak{m}_i\cdot K[N]_T$.

Now, for a semilocal Noetherian ring $R$ with Jacobson ideal $\mathfrak{m}=\bigcap_{i=1}^r \mathfrak{m}_i$ containing a field of positive characteristic, we can define the Hilbert-Kunz function as before.
For a finite $R$-module $M$ and an $\mathfrak{m}$-primary ideal $\mathfrak{n}$ we set
\[ \hkf^R(\mathfrak{n},M,q)=\len(M/\mathfrak{n}^{[q]}M). \]
If $J$ is an ideal in a Noetherian ring $R$ with $V(J)=\{\mathfrak{m}_1,\dots,\mathfrak{m}_r\}$, then
\[ \len^R(M/JM)=\len^{R_T}(M_T/JM_T) \]
for $T=\bigcap_{i=1}^r R\setminus \mathfrak{m}_i$. In this way we consider $K[N_+]$-primary ideals in $K[N]$ for a semipositive binoid $N$,
and we write $\hkf^{K[N]}(K[\mathfrak{n}],K[S],q)$ instead of $\hkf^{K[N]_T}(K[\mathfrak{n}]\cdot K[N]_T,K[S]_T,q)$.

\tr
\label{HK binoid ring}
Let $K$ be a field of characteristic $p$, $N$ a finitely generated, semipositive binoid, $S$ an $N$-set, $\mathfrak{n}$ an $N_+$-primary ideal and $q=p^e$. Then we have 
\[\hkf^N(\mathfrak{n},S,q)=\hkf^{K[N]}(K[\mathfrak{n}],K[S],q).\]
\tre
\pr
We know, from Proposition \ref{K algebra of quotient}, that
\[K[S/(S+[q]\mathfrak{n})]\cong K[S]/(K[S]\cdot K[\mathfrak{n}]^{[q]}),\] and by a dimension count we can conclude that 
\[\# (S/(S+[q]\mathfrak{n}))=\dim_K K[S/(S+[q]\mathfrak{n})]=\dim_K K[S]/(K[S]\cdot K[\mathfrak{n}]^{[q]}).\]
\pre

\tr
\label{HKM binoid ring and binoid}
Let $K$ be a field of characteristic $p$, $N$ a finitely generated, semipositive binoid. Suppose that $\dim N=\dim K[N]$ and $e_{HK}(\mathfrak{n},N)$ exists. Then 
\[e_{HK}^{K[N]}(K[\mathfrak{n}],K[N])=e_{HK}(\mathfrak{n},N)\] 
 and it is independent of the (positive) characteristic of $K$.
\tre
\pr
If we take $S=N$ in Theorem \ref{HK binoid ring} then we have 
\[\hkf^{K[N]}(K[\mathfrak{n}],K[N],q)=\hkf^N(\mathfrak{n},N,q).\] 
By assumption $(q=p^e,d=\dim N)$ 
\[e_{HK}(\mathfrak{n},N)=\lim_{q\rightarrow \infty} \!\dfrac{\hkf^N(\mathfrak{n},N,q)}{q^d}=\lim_{e\rightarrow \infty} \! \dfrac{\hkf^{K[N]}(K[\mathfrak{n}],K[N],q)}{q^d}=e_{HK}^{K[N]}(K[\mathfrak{n}],K[N]).\]
\pre

\tr[Miller conjecture for cancellative binoid rings]
\label{Miller conjecture}
Let $K$ be a field of characteristic $p$, $N$ a binoid.
Suppose that $\dim N=\dim K[N]$ and $e_{HK}(\mathfrak{n},N)$ exists. Then
 \[\lim_{p\rightarrow\infty} \frac{\hkf^{(\z/p)[N]}((\z/p)[\mathfrak{n}],(\z/p)[S],p)}{p^{\dim N}}\] exists and equals 
 \[\lim_{p\rightarrow\infty} e_{HK}^{(\z/p)[N]}((\z/p)[\mathfrak{n}],(\z/p)[S]).\]
\tre
\pr
This follows from the identity $\hkf^{(\z/p)[N]}((\z/p)[\mathfrak{n}],(\z/p)[S],p)=\hkf^N(\mathfrak{n},S,p)$ from Theorem \ref{HK binoid ring} and the existence of the limit over all numbers.
\pre

\tr
\label{existsrational}
Let $K$ be a field of characteristic $p$, $N$ be a finitely generated, semipositive, cancellative, reduced binoid and $\mathfrak{n}$ be an $N_+$-primary ideal of $N$. Then 
\[e_{HK}^{K[N]}(K[\mathfrak{n}],K[N])\] 
exists, is independent of the characteristic of $K$ and it is a rational number.
\tre
\pr
By Theorem \ref{f.g,s.p,c,r binoid} we know that $e_{HK}(\mathfrak{n},N)$ exists and that it is a rational number. So by Theorem \ref{HKM binoid ring and binoid} we have the result.
\pre



\begin{thebibliography} {10}
\addcontentsline{toc}{chapter}{Bibliography}



\bibitem [Bat14]{Batsukhthesis} B. Batsukh: \emph{Hilbert-Kunz theory for binoids}. Thesis, Osnabrück, 2014.




\bibitem [Boet15]{Simone} S. Boettger: \emph{Monoids with absorbing elements and their
associated algebras}. Thesis, Osnabrück, 2015.

\bibitem [Bre06]{Bre06} H. Brenner: \emph{The rationality of the Hilbert-Kunz multiplicity in graded
dimension two}. Math. Ann. 334, (2006), 91-110.

\bibitem [Bre07]{Bre07} Brenner, Holger: \emph{The Hilbert-Kunz function in graded dimension two.}
Comm. Algebra 35 (2007), no. 10, 3199–3213.

\bibitem [BreLiMil12]{BreLiMil} Brenner, Holger; Li, Jinjia; Miller, Claudia: 
\emph{A direct limit for limit Hilbert-Kunz multiplicity for smooth projective curves.}
J. Algebra 372 (2012), 488–504.

\bibitem [Bre13]{Brenner} H. Brenner: \emph{Irrational Hilbert-Kunz multiplicities}.
arXiv:1305.5873, 2013.

\bibitem [Bru05]{Bruns} W. Bruns: \emph{Conic divisor classes over a normal monoid algebra}.
Contemporary Mathematics. 390 (2005) 63-71.

\bibitem [Con96]{Conca} A. Conca: \emph{Hilbert-Kunz function of monomial ideals and binomial
hypersurfaces}. Manuscripta Math. 90 (1996), 287-300.

\bibitem[Eto02]{Eto} K. Eto: \emph{Multiplicity and Hilbert-Kunz multiplicity of monoid rings}.
Tokyo J.Math. 25 (2002), 241-245.

\bibitem [Kun69]{Kunz} E. Kunz: 
\emph{Characterizations of regular local rings of characteristic $p$}. Amer.
J. Math. 41(1969), 772-784.

\bibitem [Mon83]{Monsky}  P. Monsky: \emph{The Hilbert-Kunz function}, Math. Ann 263(1983), 43-49.

\bibitem [MonHan93]{MonHan} P. Monsky and C. Han. 
\emph{Some surprising Hilbert-Kunz functions.} Math. Z. 214. (1993), 119-135.

\bibitem [MonTei04]{MonTei04} Monsky, Paul; Teixeira, Pedro: 
\emph{$p$-fractals and power series. I: Some 2 variable results.}
J. Algebra 280, No. 2, 505-536 (2004).

\bibitem [MonTei06]{MonTei06} Monsky, Paul; Teixeira, Pedro: 
\emph{$p$-fractals and power series. II: Some applications to Hilbert-Kunz theory.}
J. Algebra 304, No. 1, 237-255 (2006).

\bibitem [Sei97]{Seibert}  G. Seibert: 
\emph{The Hilbert-Kunz function of rings of finite Cohen-Macaulay type},
Arch. Math. (Basel) 69(1997), 286-296.

\bibitem [TriFak03]{TriFak03} Fakhruddin, N. Trivedi, V.: 
\emph{Hilbert--Kunz functions and multiplicities for full flag varieties and elliptic curves.}
J. Pure Appl. Algebra 181, No. 1, 23-52 (2003).

\bibitem [Tri05a]{Trivedi1} Trivedi, V.: 
\emph{Strong semistability and Hilbert-Kunz multiplicity for singular plane curves.}
Ghorpade, Sudhir (ed.) et al., 
Commutative algebra and algebraic geometry. 
Joint international meeting of the American Mathematical Society and the Indian Mathematical Society, Bangalore, India, December 17--20, 2003. 
Providence, RI: American Mathematical Society (AMS). Contemporary Mathematics 390, 165-173 (2005).

\bibitem [Tri05b]{Trivedi2} V. Trivedi: 
\emph{Semistability and Hilbert-Kunz multiplicities for curves}.
J. Algebra 284 (2005), 627–644.

\bibitem [Tri07]{Trivedi07} V. Trivedi: \emph{Hilbert-Kunz multiplicity and reduction mod $p$}.
Nagoya Math. J. 185, 123-141 (2007).

\bibitem [Wat00]{Watanabe} K.-I. Watanabe: \emph{Hilbert-Kunz multiplicity of toric rings}. Proc. Inst.
Nat. Sci. 35 (2000), 173-177.

\bibitem [WatYos00]{Watanabeyoshida} K.-I. Watanabe: \emph{Hilbert-Kunz multiplicity and an inequality between	multiplicity and colength}.
J. Algebra 230 (2000), 295-317.











\end{thebibliography}
\end{document}